\documentstyle{article}

\newcommand{\mb}{\makebox}
\newcommand{\N}{\mb{$I\!\!N$}}
\newcommand{\R}{\mb{$I\!\!R$}}

\newcommand{\beqa}{\begin{eqnarray*}}
\newcommand{\eeqa}{\end{eqnarray*}}

\newcommand{\es}{\emptyset}
\newcommand{\ci}{\subseteq}

\renewcommand{\u}{\cup}
\renewcommand{\i}{\cap}

\newcommand{\lmy}[1]{\lim_{#1 \rightarrow \infty}}

\newcommand{\ra}{\rightarrow}
\newcommand{\Ra}{\Rightarrow}
\newcommand{\Lra}{\Leftrightarrow}

\renewcommand{\a}{\alpha}

\newcommand{\e}{\varepsilon}
\newcommand{\k}{\kappa}
\renewcommand{\l}{\lambda}
\newcommand{\m}{\mu}
\newcommand{\Si}{\Sigma}
\renewcommand{\c}{\chi}
\newcommand{\Om}{\Omega}

\newcommand{\pf}{\noindent {\bf Proof :} }
\newcommand{\f}{\frac}
\newcommand{\ov}{\overline}

\newcommand{\qed}{\hfill\rule{3mm}{3mm}}

\newcounter{cnt1}
\newcounter{cnt2}
\newcounter{cnt3}
\newcounter{cnt4}
\newcommand{\blr}{\begin{list}{$($\roman{cnt1}$)$} {\usecounter{cnt1}
		\setlength{\topsep}{0pt} \setlength{\itemsep}{0pt}}}
\newcommand{\blR}{\begin{list}{\Roman{cnt4}.\ } {\usecounter{cnt4}
		\setlength{\topsep}{0pt} \setlength{\itemsep}{0pt}}}
\newcommand{\bla}{\begin{list}{$($\alph{cnt2}$)$} {\usecounter{cnt2}
		\setlength{\topsep}{0pt} \setlength{\itemsep}{0pt}}}
\newcommand{\bln}{\begin{list}{$($\arabic{cnt3}$)$} {\usecounter{cnt3}
                \setlength{\topsep}{0pt} \setlength{\itemsep}{0pt}}}
\newcommand{\el}{\end{list}}

\newtheorem{ex}{Example}[section]
\newtheorem{Q}{Question}[section]

\newtheorem{Thm}{Theorem}[section]
\newtheorem{Defn}{Definition}[section]

\newtheorem{Cor}[Thm]{Corollary}

\newtheorem{rem}{Remark}[section]
\newcommand{\Rem}{\begin{rem} \rm}
\title{\bf On A New Asymptotic Norming Property}
\author{\bf Pradipta Bandyopadhyay \& Sudeshna Basu\\
\sc Stat--Math Division, Indian Statistical Institute\\ \sc 203,
B.\ T.\ Road, Calcutta 700 035, India\\ e-mail~:
pradipta@isical.ernet.in \& res9415@isical.ernet.in}
\date{}
\begin{document}

\vfill
\maketitle
\thispagestyle{empty}
\vfill

\begin{abstract}
In this work, we introduce a new Asymptotic Norming Property
(ANP) which lies between the strongest and weakest of the
existing ones, and obtain isometric characterisations of it.
The corresponding w*-ANP turns out to be equivalent on the one
hand, to Property $(V)$ introduced by Sullivan, and to a ball
separation property on the other. We also study stability
properties of this new ANP and its w*-version.
\end{abstract}
\vfill

\noindent \hrulefill\
\begin{description}
\item[AMS Subject Classification (1990) :] 46B20, 46B22
\item[Keywords and Phrases :] (w*-) Asymptotic Norming Property,
Hahn-Banach Smooth, Very Smooth, Kadec-Klee Property, Property
(V), Ball separation Properties.
\item[ISI Tech.\ Report No.] 5/95 \hfill Date~: 9/5/95
\end{description}
\pagebreak

\section{Introduction}

The Asymptotic Norming Property (ANP) was introduced by James
and Ho \cite{JH} to show that the class of separable Banach
spaces with Radon-Nikod\'ym Property (RNP) is larger than those
isomorphic to subspaces of separable duals. Three different
Asymptotic Norming Properties were introduced and were shown to
be equivalent in separable Banach spaces. Recently, Hu and Lin
\cite{HL1} have obtained isometric characterisations of the ANPs
and shown that they are equivalent in Banach spaces admitting a
locally uniformly convex renorming, a class larger than
separable Banach spaces. In dual Banach spaces, they introduced
a stronger notion called the w*-ANP, which turned out to be nice
geometric properties.

Here, we introduce a new Asymptotic Norming Property which lies
between the strongest and weakest of the existing ones, and
obtain isometric characterisations of it. The corresponding
w*-ANP, the main object of our study, turns out to be equivalent
on the one hand, to Property $(V)$ introduced by Sullivan
\cite{S}, and to a ball separation property {\em \`a la}\/ Chen
and Lin \cite{CL} on the other. We also study stability
properties of this new ANP and its w*-version.

\section{The ANPs~: Old and New}

\begin{Defn} \rm 
For a Banach space $X$, let $S_{X} = \{x~: \|x\| =1\}$ and $B_X
= \{x~: \|x\| \leq 1\}$. 

\noindent A subset $\Phi$ of $B_{X^*}$ is called a norming set
for $X$ if $\|x\| = \sup_{x^* \in \Phi} x^*(x)$, for all $x \in
X$. A sequence $\{x_n\}$ in $S_X$ is said to be asymptotically
normed by $\Phi$ if for any $\e > 0$, there exists a $x^* \in
\Phi$ and $N \in \N$ such that $x^*(x_n) > 1 - \e$ for all $n
\geq N$.

\noindent For $\k =$ I, II or III, a sequence $\{x_n\}$ in $X$
is said to have the property $\k$ if
\blR
\item $\{x_n\}$ is convergent.
\item $\{x_n\}$ has a convergent subsequence.
\item $\bigcap_{n = 1}^{\infty} \ov{co} \{x_k~:k\geq n\} \neq
\es$, where $\ov{co}(A)$ is the closed convex hull of $A \ci X$.
\el
For $\k =$ I, II or III, $X$ is said to have the asymptotic
norming property $\k$ with respect to $\Phi$ ($\Phi$-ANP-$\k$),
if every sequence in $S_X$ that is asymptotically normed by
$\Phi$ has property $\k$.
\end{Defn}

\Rem In \cite[Theorem 2.3]{HL1}, it is shown that $\Phi$-ANP-III
is equivalent to the apparently stronger property that every
sequence in $S_X$ asymptotically normed by $\Phi$ has a weakly
convergent subsequence. 
\end{rem}

This motivates the following definition.
\begin{Defn} \rm 
Let $X$ be a Banach space and let $\Phi \ci B_{X^*}$ be a
norming set for $X$. $X$ is said to have $\Phi$-ANP-II$'$ if any
sequence $\{x_n\}$ in $S_X$ which is asymptotically normed by
$\Phi$ is weakly convergent. 
\end{Defn}

\begin{Defn} \rm
$X$ is said to have the asymptotic norming property $\k$
(ANP-$\k$), $\k =$ I, II, II$'$ or III, if there exists an
equivalent norm $\|\cdot\|$ on $X$ and a norming set $\Phi$ for
($X, \|\cdot\|$) such that $X$ has $\Phi$-ANP-$\k$.
\end{Defn}

\Rem
Clearly, $\Phi$-ANP-I $\Ra$ $\Phi$-ANP-II$'$ $\Ra$
$\Phi$-ANP-III. Thus all the ANPs are equivalent in Banach
spaces admitting a locally uniformly convex renorming, in
particular, in separable Banach spaces.
\end{rem}

\begin{Defn} \rm
A Banach space $X$ is said to have the Kadec Property (K) if the
weak and the norm topologies coincide on the unit sphere, i.e.,
$(S_X, w) = (S_X, \|\cdot\|)$. 

\noindent $X$ is said to have Kadec-Klee Property (KK) if for
any sequence $\{x_n\}$ and $x$ in $B_X$ with $\lim_n \|x_n\| =
\|x\| = 1$ and w-$\lim_n x_n = x$, $\lim_n \|x_n - x\| = 0$.
\end{Defn}

The proofs of the following two theorems are evidently similar
to those of Theorems 2.4 and 2.5 of \cite{HL1}. We include the
details only when we feel some elaboration is needed.

\begin{Thm} \label{I'}
Let $\Phi$ be a norming set for a Banach space $X$. The
following are equivalent~:
\bla
\item $X$ has $\Phi$-ANP-I
\item $X$ has $\Phi$-ANP-II$'$ and $X$ has (K)
\item $X$ has $\Phi$-ANP-II$'$ and $X$ has (KK)
\el
\end {Thm}

\pf Since $\Phi$-ANP-II implies (K) \cite[Theorem 2.4, $(1) \Ra
(2)$]{HL1}, so does $\Phi$-ANP-I. Thus $(a) \Ra (b)$ follows and
$(b) \Ra (c)$ is obvious.

$(c) \Ra (a)$ Since $X$ has $\Phi$-ANP-II$'$, any sequence
$\{x_n\}$ in $S_X$ asymptotically normed by $\Phi$ is weakly
convergent to some $x\in X$. Then by \cite[Lemma 2.2]{HL1},
$\|x\| = 1$ and hence by (KK) we have $x_n \ra x$ in norm. \qed

\begin{Thm} \label{II'}
Let $\Phi$ be a norming set for a Banach space $X$. The
following are equivalent~:
\bla
\item $X$ has $\Phi$-ANP-II$'$.
\item $X$ has $\Phi$-ANP-III and $X$ is strictly convex.
\el
\end{Thm}

\pf $(a) \Ra (b)$ Strict convexity of $X$ follows similarly as
in the proof of \cite[Theorem 2.5, $(1) \Ra (2)$]{HL1}.

$(b) \Ra (a)$ Let $\{x_n\}$ be a sequence in $S_X$
asymptotically normed by $\Phi$. Since $X$ has $\Phi$-ANP-III,
$D = \bigcap\ov{co}\{x_k~: k\geq n \} \neq \es$. Now $X$ has
$\Phi$-ANP-III implies $\{x_n\}$ has weak cluster points and all
of them must be in $D$. Since $D\ci S_X$ and $X$ is strictly
convex, $D$ is a singleton. Moreover, since every subsequence of
$\{x_n\}$ is also asymptotically normed by $\Phi$, that
singleton is the weak limit of $\{x_n\}$. Hence $X$ has
$\Phi$-ANP-II$'$.
\qed

Some renorming results similar to Theorem 2.7 of \cite{HL1}
can easily be obtained from our results. But in this work, we
concentrate on the ANPs as isometric properties.

\section{w*-ANPs}

\begin{Defn} \rm
Let $X^*$ be a dual Banach space. $X^*$ is said to have
w*-ANP-$\k$ ($\k =$ I, II, II$'$ or III) if there exists an
equivalent norm $\|\cdot\|$ on $X$ and a norming set $\Phi$ for
$X^*$ in $B_X$ such that $X^*$ has $\Phi$-ANP-$\k$.
\end{Defn}

\Rem If $\Phi \ci B_X$ is a norming set for $X^*$, then
$\ov{co}(\Phi \u - \Phi) = B_X$. Hence, by \cite[Lemma 3]{HL3}
and similar arguments, $\Phi$-ANP-$\k$ is equivalent to
$B_X$-ANP-$\k$ ($\k =$ I, II, II$'$ or III). Thus, we can and do
work with $\Phi = B_X$.
\end{rem}

\begin{Defn} \rm
For a Banach space $X$, let $X^{\perp} = \{x^{\perp} \in
X^{***}$~: $x^{\perp}(x) = 0$ for all $x\in X\}$. A Banach space
is said to be Hahn-Banach smooth if for all $x^*\in X^*$, $\|x^*
+ x^{\perp}\| = \|x^*\| = 1$ implies $x^{\perp} = 0$, i.e., $x^*
\in X^{***}$ is the unique norm preserving extension of
$x^*|_X$.
\end{Defn}

\Rem \label{r1} It is shown in \cite[Theorem 1]{HL2} that $X$ is
Hahn-Banach smooth if and only if $X^*$ has $B_X$-ANP-III if and
only if $(S_{X^*}, w^*) = (S_{X^*}, w)$.
\end{rem}

Extending a characterisation of rotundity of $X^*$ due to L.\
P.\ Vlasov \cite{V}, Sullivan \cite{S} introduced the following
stronger property~:
\begin{Defn} \rm
A Banach space $X$ is said to have the Property $(V)$, if there
do not exist an increasing sequence $\{B_n\}$ of open balls with
radii increasing and unbounded, and norm one functionals $x^*$
and $y^*_k$ such that for some constant $c$,

$x^*(b) > c$ for all $b\in \u{B_n}$,

$y^*_k(b) > c$ for all $b \in B_n$, $n \leq k$ and

dist($co(y^*_1, y^*_2, \ldots), x^*) > 0$.
\end{Defn}

\begin{Defn} \rm 
Let $W \ci X^*$ be a closed bounded convex set.
\bla
\item A point $x^* \in W$ is said to be a weak*-weak point of
continuity (w*-w pc) of $W$ if $x^*$ is a point of continuity of
the identity map from $(W, w^*)$ to $(W, w)$
\item A point $x^* \in W$ is said to be a w*-strongly extreme
point of $W$ if the family of w*-slices containing $x^*$ forms a
local base for the weak topology of $X^*$ at $x^*$ (relative to
$W$).
\el
\end{Defn}

Now we have our main characterisation theorem.
\begin{Thm} \label{w*-II'}
For a Banach space $X$, the following are equivalent~:
\bla
\item $X^*$ has $B_X$-ANP-II$'$.
\item $X^*$ is strictly convex and $X$ is Hahn-Banach smooth.
\item $X$ has Property $(V)$.
\item All points of $S_{X^*}$ are w*-strongly extreme points
of $B_{X^*}$.
\el
\end{Thm}

\pf $(a) \Lra (b)$ is immediate from Theorem~\ref{II'} and
Remark~\ref{r1}, while $(b) \Lra (c)$ is just \cite[Theorem
4]{S}.

$(b) \Ra (d)$ Since $(S_{X^*}, w^*) = (S_{X^*}, w)$, and the
norm is lower semi-continuous with respect to both weak and
weak* topology of $X^*$, any $x^* \in S_{X^*}$ is a w*-w pc of
$B_{X^*}$. Now, since $X^*$ is strictly convex, every $x^* \in
S_{X^*}$ is an extreme point of $B_{X^*}$. By a classical result
of G.\ Choquet \cite[Proposition 25.13]{C}, for any $x^* \in
S_{X^*}$, the family of w*-slices containing $x^*$ forms a local
base for the w*, and therefore the weak, topology of $X^*$
relative to $B_{X^*}$.

$(d) \Ra (b)$ From (d), it is immediate that $X^*$ is strictly
convex and any $x^* \in S_{X^*}$ is a w*-w pc of $B_{X^*}$. \qed

\begin{Defn} \rm
(a) The duality mapping $D$ for a Banach space $X$ is the
set-valued map from $S_X$ to $S_{X^*}$ defined by
\[D(x) = \{x^* \in S_{X^*}~: x^* (x) = 1\}, \quad x\in S_X.\]

(b) A Banach space $X$ is said to be very smooth if every $x\in
S_X$ has a unique norming element in $X^{***}$.
\end{Defn}

It is known that $X$ is Fr\'echet differentiable (very smooth)
if and only if the duality mapping D is single-valued and is
norm-norm (norm-weak) continuous.

\begin{Thm} \label{VS}
If $X^*$ has $B_X$-ANP-II$'$, then $X$ is very smooth.
\end{Thm}

\pf That Property $(V)$ implies very smooth was already observed
in \cite{S}. We, however, prefer the following direct and
ANP-like argument similar to \cite[Theorem 4(1)]{HL2}.

Since $X^*$ is strictly convex, $X$ is smooth. Now let $\{x_n\}
\ci S_X$ be such that $x_n \ra x$. Let $\{x^*_n\} = D(x_n)$, we
have $|x^*_n (x) - 1| \leq |x^*_n(x) - x^*_n(x_n)| \leq
\|x^*_n\| \|x - x_n\| \leq \|x - x_n \| \ra 0$ as $n\ra \infty$.
That is, $\lmy {n} x^*_n(x) = 1$. So $\{x^*_n\} $ is
asymptotically normed by $B_X$, and hence, is weakly convergent
to $x^*$ (say). Clearly, $x^* \in D(x)$ and since $X$ is smooth,
$\{x^*\} = D(x)$. Hence $X$ is very smooth.
\qed

Analogous to \cite{HL2}, we have the following question~:
\begin{Q} \label{q2}
Let $X$ be Hahn-Banach smooth and very smooth. Does $X^*$ have
$B_X$-ANP-II$'$?
\end{Q}

\begin{ex}
In general, for a Banach space $X$, the properties ANP-I, II,
II$'$ and III are all distinct, i.e., except for the obvious
implications $\Phi$-ANP-I $\Ra \Phi$-ANP-II $\Ra \Phi$-ANP-III
and $\Phi$-ANP-I $\Ra \Phi$-ANP-II$' \Ra \Phi$-ANP-III, no other
implication is generally true.
\end{ex}

\pf Clearly, it suffices to show that none of $\Phi$-ANP-II or
$\Phi$-ANP-II$'$ implies the other.

(1) Let $X = c_0$, $X^* = \ell_1$. Since $(S_{X^*}, w) =
(S_{X^*}, \|\cdot\|)$ on $\ell_1$, by \cite[Theorem 3.1]{HL1},
$X^*$ has $B_X$-ANP-II. But $X^*$ is not strictly convex.

(2) \cite{S} On $X = \ell_2$, define an equivalent norm as
$\|x\|_0 = \max\{1/2(\|x\|_2)$, $\|x\|_{\infty}\}$. And define
$T~: \ell_2 \ra \ell_2$ by $T(\a_k) = \a_k/k$, for $(\a_k)\in
\ell_2$. Then $T$ is an 1-1 continuous linear map. Hence the
equivalent dual norm $\|x\|_3 = \|x\|_0 + \|Tx\|_2$ is strictly
convex \cite{D1}. Also since $\ell_2$ is reflexive, it has
$B_X$-ANP-III with respect to $\|\cdot\|_3$. Thus by
Theorem~\ref{II'}, $(\ell_2, \|\cdot\|_3)$ has $B_X$-ANP-II$'$.
But, it was observed in \cite{S} that $(\ell_2, \|\cdot\|_3)$
lacks (KK). \qed

\begin{ex}
The above two examples show that a space may have ANP-III, but
may lack either ANP-II or II$'$. The following is an example of
a Banach space which has ANP-III but lacks both ANP-II and
II$'$.
\end{ex}

\pf Let $X = \ell_2 \oplus_1 \R$. It is clear that $X^* = \ell_2
\oplus_{\infty} \R$ is reflexive, and hence, has $B_X$-ANP-III.
However $X^*$ is not strictly convex, and hence cannot have
$B_X$-ANP-II$'$. Also the weak and the norm topologies do not
coincide on $S_{X^*}$. Indeed, since $\ell_2$ is infinite
dimensional, by Riesz' Lemma, there exists a sequence $\{x_n\}$
in $S_{\ell_2}$ such that $\|x_n - x_m\|_2 \geq 1$, $n \neq m$.
Let $z_n = (x_n, 1)$. So $\|z_n\|_{\infty} = 1$ and $\|z_n -
z_m\|_{\infty} \geq 1$. Clearly, $\{z_n\}$ cannot have any norm
convergent subsequence. But as $\ell_2$ is reflexive, $\{x_n\}$
has a subsequence $\{x_{n_k}\}$ converging weakly to some $x \in
B_{\ell_2}$. Then obviously $(x_{n_k}, 1) = z_{n_k}$ converges
weakly to $(x, 1) = z$ (say) and $\|z\|_{\infty} = 1$.
\qed

\section{A Ball Separation Property}

In a recent work, Chen and Lin \cite{CL} have obtained certain
ball separation properties which, in equivalent formulations,
characterise $B_X$-ANP-$\k$ ($\k =$ I, II and III).

Here we obtain a similar characterisation of $B_X$-ANP-II$'$. In
fact, as in \cite{CL}, we also take a local approach, i.e., we
characterise w*-strongly extreme points of $B_{X^*}$. And for
that we need the following characterisation for w*-w pc of
$B_{X^*}$ which is immediate from Theorem 3.1 of \cite{CL}.
\begin{Thm} \label{pc}
For a Banach space $X$ and $f_0 \in S_{X^*}$, the following are
equivalent~:
\blr
\item $f_0$ is a w*-w pc of $B_{X^*}$.
\item for any $x_0^{**} \in X^{**}$ and $\a \in \R$, if
$f_o(x_0^{**}) > \a$, then there exists a ball $B^{**}$ in
$X^{**}$ with centre in $X$ such that $x_0^{**} \in B^{**}$ and
$\inf f_0(B^{**}) > \a$.
\el
\end{Thm}

>From \cite[Theorem 1.3]{CL} and the arguments of \cite[Corollary
2]{B}, we get
\begin{Thm} \label{Ext}
For a Banach space $X$ and $f_0 \in S_{X^*}$, the following are
equivalent~:
\blr
\item $f_0$ is an extreme point of $B_{X^*}$.
\item for any compact set $A \ci X$ if $\inf f_0(A) > 0$, then
there exists a ball $B$ in $X$ such that $A \ci B$ and $\inf
f_0(B) > 0$.
\item for any finite set $A \ci X$ if $\inf f_0(A) > 0$, then
there exists a ball $B$ in $X$ such that $A \ci B$ and $\inf
f_0(B) > 0$.
\el
\end{Thm}

\begin{Thm}
For a Banach space $X$ and $f_0 \in S_{X^*}$, the following are
equivalent~:
\bla
\item $f_0$ is a w*-strongly extreme point of $B_{X^*}$.
\item $f_0$ is a w*-w pc and an extreme point of $B_{X^*}$.
\item for any compact set $A \ci X^{**}$, if $\inf f_0(A) > 0$,
then there exists a ball $B^{**} \ci X^{**}$ with centre in $X$
such that $A \ci B^{**}$ and $\inf f_0(B^{**}) > 0$.
\item for any finite set $A \ci X^{**}$, if $\inf f_0(A) > 0$,
then there exists a ball $B^{**} \ci X^{**}$ with centre in $X$
such that $A \ci B^{**}$ and $\inf f_0(B^{**}) > 0$.
\el
\end{Thm}

\pf $(a) \Lra (b)$ is just the local version of
Theorem~\ref{w*-II'} $(b) \Lra (d)$.

$(b) \Ra (c)$ Since $f_0$ is a w*-strongly extreme point of
$B_{X^*}$, it is easily seen that it remains extreme in
$B_{X^{***}}$. Thus by Theorem~\ref{Ext}, for any compact set
$A$ in $X^{**}$ with $\inf f_0(A) > 0$, there exists a ball in
$B^{**} = B^{**}(x^{**}, r) \ci X^{**}$ such that $A \ci B^{**}$
and $\inf f_0(B^{**}) > 0$. Now, $\inf f_0(B^{**}(x^{**}, r)) >
0$ implies $f_0(x_0^{**}) > r$. Since $f_0$ is a w*-w pc, by
Theorem~\ref{pc}, there exists a ball $B^{**}(x,t) \ci X^{**}$
such that $x_0^{**} \in B^{**}(x,t)$ and $\inf f_0(B^{**}(x,t) >
r$. This implies $f_0(x) > r + t$. Thus, $A \ci B^{**}(x_0^{**},
r) \ci B^{**}(x, r + t)$ and $\inf f_0(B^{**}(x, r + t) > 0$.

$(c) \Ra (d)$ is trivial.

$(d) \Ra (b)$ Taking $A \ci X$, it follows from
Theorem~\ref{Ext} that $f_0$ is extreme in $B_{X^*}$. And taking
$A$ to be a singleton, it follows from Theorem~\ref{pc} that
$f_0$ is an w*-w pc.
\qed

\begin{Cor}
For a Banach space $X$, the following are equivalent~:
\blr
\item $X^*$ has $B_X$-ANP-II$'$.
\item for any w*-closed hyperplane $H$ in $X^{**}$, and any
compact convex set $A$ in $X^{**}$ with $A \i H = \es$, there
exists a ball $B^{**}$ in $X^{**}$ with centre in $X$ such that
$A \ci B^{**}$ and $B^{**} \i H = \es$.
\item for any w*-closed hyperplane $H$ in $X^{**}$, and any
finite dimensional convex set $A$ in $X^{**}$ with $A \i H =
\es$, there exists a ball $B^{**}$ in $X^{**}$ with centre in
$X$ such that $A \ci B^{**}$ and $B^{**} \i H = \es$.
\el
\end{Cor}

\section{Stability Results}

\begin{Thm} 
Let $X$ be a Banach space with $\Phi$-ANP-$\k$, $\k =$ I, II,
II$'$ or III. Then any closed subspace $Y$ of $X$ has
$\Phi|_Y$-ANP-$\k$, where $\Phi|_Y = \{y^*~: y^* = x^*|_Y$, $x^*
\in \Phi\}$.
\end{Thm}

\begin{Thm} 
Let $X$ be a Banach space such that $X^*$ has
$B_X$-ANP-$\k$, $\k =$ I, II, II$'$ or III. Then for any
closed subspace $Y$ of $X$, $Y^*$ has $B_Y$-ANP-$\k$.
\end{Thm}

\pf Let $\{y^*_n\} \ci S_{Y^*}$ be asymptotically normed by
$B_Y$. For every $n \geq 1$, let $x^*_n$ be a norm preserving
extension of $y^*_n$ to $X$. Then $\{x^*_n\}$ is asymptotically
normed by $B_X$, and hence has property $\k$. Now the
restriction map $x^* \ra x^*|_Y$ brings property $\k$ back
to $\{y^*_n\}$. \qed

\begin{Cor} \label{her}
Hahn-Banach smoothness and Property $(V)$ are hereditary.
\end{Cor}

\Rem
This observation appears to be new. Note that we do not need the
stability of the ANPs under quotients to prove the above
theorem. In fact, it is not clear whether the ANPs are indeed
stable under quotients.
\end{rem}

Let $X$ be a Banach space, $1< p, q < \infty$ with $ 1/p + 1/q =
1$ and $(\Om ,\Si, \m)$ be a positive measure space so that
$\Si$ contains an element with finite positive measure. Let
$\Phi$ be a norming set for $X$. Then define $\Phi_1 = co(\Phi
\u \{0\}) \setminus S_X$ and
\beqa 
\Delta _n & = & \{\sum_{i=1}^{m} \l_i x ^*_i \c_{E_i}~: x^*_i
\in \Phi_1, \f{(n-1)}{n} \leq \|x^*_i\| \leq \f{n}{(n+1)}, E_i
\in \Si, \\ && E_i \i E_j = \es, \mb{ for } i \neq j, \l_i > 0
\mb{ with } \sum_{i =1}^{m} \l^q_i \m(E_i) = 1\}
\eeqa
Then $\Delta ( \Phi, \m, q) = \bigcup_{n \geq 1} \Delta_n$ is a
norming set for $L^p ( \m, X)$ \cite{HL3}.

\begin{Thm}
Let $X$ be a Banach space, $\Phi \ci B_{X^*}$ be a norming set
for $X$. $X$ has $\Phi$-ANP-II$'$ if and only if $L^p (\m, X)$
has $\Delta( \Phi, \m, q)$-ANP-II$'$.
\end{Thm}

\pf It is well-known that $X$ is strictly convex if and only if
$L^p (\m, X)$ is strictly convex \cite{D2}. And in \cite[Theorem
6]{HL3}, it is shown that $X$ has $\Phi$-ANP-III if and only if
$L^p (\m, X)$ has $\Delta( \Phi, \m, q)$-ANP-III. Now, the
result follows from Theorem~\ref{II'}. \qed

\Rem 
Let $X$ be a Banach space. If $(X, \|\cdot\|)$ has ANP-II, the
space $(L_p(\m, X), \|\cdot\|)$ may not have ANP-II. For an
example, see \cite{HL3}. Thus we have nicer stability results
for ANP-II$'$ which was lacking in ANP-II.
\end{rem}

\begin{Thm}
Let $X$ be a Banach space. $X$ has Property $(V)$ if and only
if $L^p(\m, X)$ has $(V)$ $(1< p < \infty)$.
\end{Thm}

\pf By Corollary~\ref{her}, $X$ inherits Property $(V)$ from
$L^p(\m, X)$.

Conversely, if $X$ has Property $(V)$, by \cite[Theorem 4]{S},
$X$ is Hahn-Banach smooth. Hence $X$ is an Asplund space. Thus,
$L^p(\m, X)^* = L^q(\m, X^*)$, where $1/p + 1/q = 1$. From
\cite[Theorem 6]{HL3}, $L^p(\m, X)$ is Hahn-Banach smooth. Also,
$X^*$ strictly convex implies $L^q(\m, X^*)$ is strictly convex.
The result now follows from Theorem~\ref{II'}. \qed

\end{document}